\documentclass[notitlepage,leqno,11pt]{article}
\usepackage{amssymb}
\catcode`\@=11
\@addtoreset{equation}{section}

\catcode`\@=12

\usepackage{latexsym}
\usepackage{amsmath}
\usepackage{amsfonts}
\usepackage{amssymb}
\usepackage{mathrsfs}

\renewcommand{\b }{\beta }

\renewcommand{\O }{\Omega }
\newcommand{\be}{\begin{equation}}
\newcommand{\ee}{\end{equation}}

\newcommand{\R}{\mathbb{R}}
\newcommand{\de}{\partial}
\def\C{{\mathcal C}}
\def\S{{\mathbb S}}
\def\f{{\varphi}}
\def\div{{\rm div}}
\def\eps{{\varepsilon}}

\def\D{\mathcal D^{1,2}}
\def\fn{\varphi_{\delta,n}}
\def\vn{v_{\delta,n}}
\def\tf{\tilde\varphi}
\newcommand{\n}{\nabla}
\newcommand{\wt }{\widetilde}
\newcommand{\Sig }{\Sigma}
\newcommand{\calC }{\mathcal{C}}

\newcommand{\weight}{|x|^{-2}\left|\log|x|\right|^{-2}}
\newcommand{\Gweight}{|x|^{-2}\left|\log|x|\right|^{-2\gamma}}
\newcommand{\CR}{{\C_\Sigma^R}}
\newcommand{\sobolev}{\widetilde H^1}

\newcommand{\hardy}{\frac{(N-2)^2}{4}}

\newtheorem{Theorem}{Theorem}[section]
\newtheorem{Lemma}[Theorem]{Lemma}
\newtheorem{Proposition}[Theorem]{Proposition}

\newtheorem{Remark}[Theorem]{Remark}

\def\proof{\noindent{{\bf Proof. }}}
\def\square{\vbox{
    \hrule height .4pt
    \hbox{\vrule width .4pt height 7pt \kern 7pt
       \vrule width .4pt}
    \hrule height .4pt }}

\def\QED{\hfill {$\square$}\goodbreak \medskip}

\begin{document}

\title{ Sharp nonexistence results for a linear elliptic inequality involving Hardy and Leray potentials}

\author{Mouhamed Moustapha Fall\footnote{\footnotesize{Universit\'e Catholique de Louvain-La-Neuve,
D\'epartement de Math\'ematique.
 Chemin du Cyclotron 2,  1348 Louvain-la-Neuve, Belgique.
E-mail: {\tt mouhamed.fall@uclouvain.be.}  } } ~and Roberta
Musina\footnote{\footnotesize{Dipartimento di Matematica ed
Informatica, Universit\`a di Udine, via delle Scienze, 206 -- 33100
Udine, Italy. E-mail: {\tt musina@dimi.uniud.it.} }}}
\date{}

\maketitle

\bigskip

\noindent {\footnotesize{\bf Abstract.}
In this paper we deal with non-negative distributional supersolutions for a class of linear
elliptic equations involving inverse-square potentials and logarithmic weights.
We prove sharp nonexistence results.}
\bigskip\bigskip

\noindent{\footnotesize{{\it Key Words:} Hardy inequality, logarithmic weights, nonexistence.}}

\noindent{\footnotesize{{\it 2000 Mathematics Subject Classification:}
35J20, 35J70, 35B33}}

\bigskip
\bigskip
\vfill\eject

\section*{Introduction}
In recent years a great deal work  has been made to find necessary
and sufficient conditions for the existence of distributional
supersolutions to semilinear elliptic equations with inverse-square
potentials. We quote for instance \cite{BDT} (and the references
therein), where a problem  related to the Hardy and Sobolev
inequalities has been studied. In the present paper we are
interested in a class of  linear elliptic equations.

Let $N\ge 2$ be an integer, $R\in(0,1]$ and let  $B_R$ be the ball
in $\R^N$ of radius $R$ centered at $0$. We focus our attention on non-negative distributional solutions
to
\begin{equation}
\label{eq:problem_ball}
-\Delta u -\frac{(N-2)^2}{4}|x|^{-2}~\! u\ge \alpha |x|^{-2}\left|\log|x|\right|^{-2}~\!u\quad
\textrm{in $\mathcal D'(B_R\setminus\{0\})$,}
\end{equation}
where $\alpha\in\R$ is a varying parameter. By a standard definition, a solution to (\ref{eq:problem_ball})
is a function $u\in L^1_{\textrm loc}(B_R\setminus\{0\})$ such that
$$
-\int_{B_R}
u\Delta\f~dx-\frac{(N-2)^2}{4}\int_{B_R}|x|^{-2}u\f~dx\ge\alpha
\int_{B_R}\weight u\f~dx
$$
for any non-negative $\f\in C^\infty_c(B_R\setminus\{0\})$.
Problem (\ref{eq:problem_ball})
is motivated by the inequality
\begin{equation}
\label{eq:HL} \int_{B_1}|\nabla u|^2~dx-\frac{(N-2)^2}{4}
\int_{B_1}|x|^{-2}|u|^2\\
\ge
\frac{1}{4}\int_{B_1}|x|^{-2}\left|\log|x|\right|^{-2}|u|^2~dx~\!,
\end{equation}
which holds for any $u\in C^\infty_c(B_1\setminus\{0\})$ (see for
example \cite{ACR},  \cite{BFT}, \cite{Ch}, \cite{GM} and Appendix
\ref{S:A}). Notice that (\ref{eq:HL}) improves the Hardy inequality
for maps supported by the unit ball if $N\ge 3$. Inequality
(\ref{eq:HL}) was firstly proved by Leray \cite{Le}  in the lower
dimensional case $N=2$.

Due to the sharpness of the constants in (\ref{eq:HL}), a necessary and sufficient
condition for the existence of non-trivial and non-negative solutions
to (\ref{eq:problem_ball}) is that $\alpha\le 1/4$ (compare with
Theorem \ref{th:1/4} in Appendix
\ref{A:AP} and with Remark \ref{R:ce}).

In case $\alpha\le 1/4$ we provide necessary conditions on
the parameter $\alpha$
to have the existence of non-trivial solutions satisfying suitable integrability properties.
\begin{Theorem}
\label{th:ne_ball}
Let $R\in(0,1]$ and let $u\ge 0$ be a distributional solution to
(\ref{eq:problem_ball}). Assume that there exists $\gamma\le 1$ such that
$$
u\in L^2_{\rm loc}(B_R;\Gweight~dx)~,\quad \alpha\ge \frac{1}{4}-(1-\gamma)^2~\!.
$$
Then $u=0$ almost everywhere in $B_R$.
\end{Theorem}
We remark that Theorem \ref{th:ne_ball}  is sharp, in view of the explicit counter-example in
Remark \ref{R:ce}.

\medskip

Let us point out some consequences of Theorem \ref{th:ne_ball}. We use the Hardy-Leray inequality (\ref{eq:HL}) to
introduce the space
$\sobolev_0(B_1)$ as the
closure of $C^\infty_c(B_1\setminus\{0\})$ with respect to the
scalar product
$$
\langle u,v\rangle=\int_{B_1}\nabla u\cdot\nabla v~dx-\hardy\int_{B_1}
|x|^{-2}uv~dx
$$
(see for example \cite{DD}). It turns out that
 $\sobolev_0(B_1)$ strictly contains the standard Sobolev space $H^1_0(B_1)$,
 unless $N=2$.

Take $\gamma=1$ in Theorem \ref{th:ne_ball}. Then problem (\ref{eq:problem_ball})
has no non-trivial and non-negative solutions $u\in L^2_{\rm
loc}(B_R;\weight~dx)$ if $\alpha= 1/4$. Therefore, if in the dual
space $\sobolev_0(B_R)'$, a function $u\in \sobolev_0(B_R)$ solves
$$
\begin{cases}
-\Delta u -\displaystyle\frac{(N-2)^2}{4}~\!|x|^{-2}~\! u\ge
 \displaystyle\frac{1}{4}~\! |x|^{-2}\left|\log|x|\right|^{-2}~\!u& \textrm{in $B_R$}\\
 u\ge 0~\!,
 \end{cases}
 $$
then  $u=0$ in $B_R$.

Next take $\gamma=0$ and $\alpha\ge -3/4$. From Theorem
\ref{th:ne_ball} it follows that  problem (\ref{eq:problem_ball})
has no non-trivial and non-negative solutions $u\in L^2_{\rm
loc}(B_R;|x|^{-2}~dx)$. In particular, if $N\geq 3$ and if $u\in
H^1_0(B_R)\hookrightarrow L^2(B_R;|x|^{-2}~dx)$ is a weak solution to
$$
\begin{cases}
-\Delta u -\displaystyle\frac{(N-2)^2}{4}~\!|x|^{-2}~\! u\ge
-\frac{3}{4}~\!\weight u &
\textrm{in $B_R$}\\
u\ge 0~\!,
\end{cases}
$$
then $u= 0$ in $B_R$.
Thus
Theorem \ref{th:ne_ball} improves some of the nonexistence results
in  \cite{AS} and in \cite{Gz}.

\medskip

The case of boundary singularities has been little studied.
In Section \ref{S:cone} we prove sharp nonexistence results for
inequalities in cone-like domains in $\R^N$, $N\ge 1$, having a
vertex at $0$. A special case concerns linear problems in
half-balls. For $R>0$ we let $B^+_R=B_R\cap \R^N_+$, where $\R^N_+$
is any half-space. Notice that $B_R^+=(0,R)$ or $B_R^+=(-R,0)$ if
$N=1$. A necessary and sufficient condition for the existence of non-negative and non-trivial
distributional solutions to
\begin{equation}
\label{eq:problem_half}
-\Delta u -\frac{N^2}{4}~\!|x|^{-2}~\! u\ge \alpha |x|^{-2}\left|\log|x|\right|^{-2}~\!u\quad
\textrm{in $\mathcal D'(B_R^+)$}
\end{equation}
is that $\alpha\le 1/4$ (see Theorem \ref{th:1/4cone} and Remark \ref{rem:sharp_cone}),
and  the following result holds.

\begin{Theorem}
\label{th:half_ball}
Let $R\in(0,1]$, $N\ge 1$ and let $u\ge 0$ be a distributional solution to (\ref{eq:problem_half}).
Assume that there exists $\gamma\le 1$ such that
$$
u\in L^2(B_R^+;\Gweight~dx)~,\quad \alpha\ge \frac{1}{4}-(1-\gamma)^2~\!.
$$
Then $u=0$ almost everywhere in $B_R^+$.
\end{Theorem}

The key step in our proofs consists in studying the
ordinary differential inequality
\begin{equation}
\label{eq:psi_intro}
\begin{cases}
-\psi''\ge \alpha s^{-2}\psi&\quad\textrm{in $\mathcal{D}'(a,\infty)$}\\
\psi\geq0~\!,
\end{cases}
\end{equation}
where $a>0$.
In our crucial Theorem \ref{th:ne_psi} we prove a nonexistence result for (\ref{eq:psi_intro}), under suitable
weighted integrability
assumptions on $\psi$. Secondly, thanks to an "averaged Emden-Fowler transform", we show
 that distributional solutions to problems of the form
(\ref{eq:problem_ball}) and (\ref{eq:problem_half}) give rise to
solutions of (\ref{eq:psi_intro}), see Section \ref{S:proof} and
\ref{S:cone} respectively. Our main existence results readily
follow from Theorem \ref{th:ne_psi}. A similar idea, but with a
different functional change, was already used in \cite{BM} to obtain
nonexistence results for a large class of superlinear problems.

\medskip

In Appendix \ref{S:A} we give a simple proof of the Hardy-Leray inequality
for maps with support in cone-like domains that includes (\ref{eq:HL}) and that motivates
our interest in problem (\ref{eq:problem_half}).

Appendix \ref{A:AP}  deals in particular with the case $\alpha>1/4$. The nonexistence
Theorems \ref{th:1/4} and \ref{th:1/4cone} follow from
an Allegretto-Piepenbrink type result (Lemma  \ref{lem:AP}).

In the last appendix we point out some related results and some
consequences of our main theorems.

\bigskip
\small
\noindent
{\bf Notation}
\\
We denote by $\R_+$ the half real line $(0,\infty)$. For $a>0$ we put
$I_a=(a,\infty)$.

We denote by
$|\O|$ the Lebesgue measure of the domain $\Omega\subset \R^N$.
Let $q\in[1,+\infty)$ and let $\omega$ be a non-negative measurable function on $\O$. The weighted
Lebesgue space $L^q(\Omega;\omega(x)~\! dx)$ is the space of measurable maps
$u$ in $\Omega$ with finite norm $\left(\int_{\Omega}|u|^q
\omega(x)~\!dx\right)^{1/q}$.  For $\omega\equiv 1$ we simply write
$L^q(\Omega)$. We embed $L^q(\Omega;\omega(x) ~\!dx)$  into $L^q(\R^N;\omega(x)
~\!dx)$ via null extension.

\normalsize

\section{Proof of Theorem \ref{th:ne_ball}}

The proof consists of two steps. In the first one we prove a nonexistence result
for a class of linear ordinary differential inequalities that might have some interest
in itself.

\subsection{Nonexistence results for problem (\ref{eq:psi_intro})}
\label{S:psi}
We start by fixing some
terminologies.
 Let $\D(\R_+)$ be the Hilbert space obtained via
 the Hardy inequality
 \begin{equation}
\label{eq:1Hardy}
\int_0^\infty|v'|^2~ds\ge \frac{1}{4}~\!\int_0^\infty s^{-2}|v|^2~ds
~,\quad~v\in C^\infty_c(\R_+)
\end{equation}
as the completion of $C^\infty_c(\R_+)$ with respect to
the scalar product
$$
\langle v,w \rangle =\int_0^\infty v'w'~ds~\!.
$$
Notice that $\D(\R_+)\hookrightarrow L^2(\R_+;s^{-2}~ds)$ with a continuous embedding and
moreover
$\D(I_a)\subset C^0(\R_+)$ by Sobolev
embedding theorem.
By H\"older inequality, the space
 $L^2(\R_+;s^2~ds)$ is continuously embedded into
the dual space $\D(\R_+)'$.

Finally, for any $a>0$ we put   $I_a=(a,\infty)$ and
$$
 \D(I_a)=\{ v\in \D(\R_+)~|~v(a)=0~\}~\!.
$$

We need two technical lemmata.

 \begin{Lemma}
 \label{L:D12}
 Let $f\in  L^2(I_a;s^{2}~ds)$ and
 $v\in C^2(\R_+)\cap L^2(I_a;s^{-2}~ds)$ be a function satisfying
 $v(a)=0$ and
\be\label{eq:vstfh}
  -v''\leq f\quad \textrm{ in } I_a~\!.
\ee Put $v^+:=\max\{v,0\}$. Then $v^+\in \D(I_a) $ and
\begin{equation}
\label{eq:vn+} \int_a^\infty|(v^+)'|^2~ds\le \int_a^\infty
fv^+~ds~\!.
\end{equation}
 \end{Lemma}

 \proof
We first show that $(v^+)'\in L^2(\R)$ and that \eqref{eq:vn+} holds. Let $\eta\in C^\infty_c(\R)$ be a
cut-off function satisfying
$$
0\le \eta\le 1~,\quad \eta(s)\equiv 1~\textrm{for $|s|\le 1$}~,\quad
\eta(s)\equiv 0~\textrm{for $s\ge 2$}
$$
and put
$\eta_h(s)=\eta(s/h)$. Then $\eta_h v^+\in \D(I_a)$ and $\eta_h v^+\ge 0$.
Multiply \eqref{eq:vstfh} by $\eta_h v^+$ and
integrate by parts to get
\begin{equation}
\label{eq:v+-1}
\int_a^{\infty}\eta_h|(v^+)'|^2~ds-\frac{1}{2}\int_a^{\infty}\eta_h''|v^+|^2~ds\leq
\int_a^{\infty}\eta_hfv^+~ds~\!.
\end{equation}
Notice that for some constant $c$ depending only on $\eta$ it results
$$
\left|\int_a^{\infty}\eta_h''|v^+|^2~ds~\!\right|\leq c\int_h^{2h} s^{-2}|v^+|^2~ds\to 0
$$
as $h\to \infty$, since $v^+\in L^2(I_a;s^{-2}~ds)$.
Moreover,
$$\displaystyle\int_a^\infty\eta_h fv^+~ds\to \int_a^\infty fv^+~ds
$$
by Lebesgue theorem, as $fv^+\in L^1(I_a)$ by H\"older inequality. In
conclusion, from (\ref{eq:v+-1}) we infer that
\begin{equation}
\label{eq:vn+0} \int_a^h|v'_+|^2~ds\le\int_a^\infty fv^+~ds+o(1)
\end{equation}
since $\eta_h\equiv 1$ on $(a,h)$. By Fatou's Lemma we get that
$(v^+)'\in L^2(I_a)$ and (\ref{eq:vn+}) readily follows from
(\ref{eq:vn+0}). To prove that $v^+\in\D(I_a)$, it is enough to
notice that $\eta_h v^+\to v^+$ in $\D(I_a)$. Indeed,
\begin{gather*}
\int_a^\infty|1-\eta_h|^2|(v^+)'|^2\le \int_h^\infty|(v^+)'|^2~ds=o(1)\\
\int_a^\infty|\eta'_h|^2|v^+|^2~ds\le c\int_h^\infty s^{-2}|v^+|^2~ds=o(1)
\end{gather*}
as $(v^+)'\in L^2(I_a)$ and $v^+\in L^2(I_a; s^{-2}~ds)$.\QED

Through the paper we let $(\rho_n)$ to be a
standard mollifier sequence in $\R$, such that the support of
$\rho_n$ is contained in the interval
$(-\frac{1}{n},\frac{1}{n})$.

\begin{Lemma}
\label{L:g} Let $a>0$ and $\psi\in L^2(I_a; s^{-2}~ds)$. Then
$\rho_n\star\psi\in L^2(I_a;s^{-2}~ds)$ and
\begin{gather}
\label{eq:g1}
\rho_n\star\psi\to \psi
\quad\textrm{in $L^2(I_a;s^{-2}~ds)$},\\
\label{eq:g2} g_n:=\rho_n\star(s^{-2}\psi)-s^{-2}(\rho_n\star\psi)
\to 0 \quad\textrm{in $L^2(I_a;s^{2}~ds)$}~\!.
\end{gather}
\end{Lemma}

\proof  We start by noticing that $\rho_n\star\psi\to \psi$ almost
everywhere. Then we use H\"older inequality to get
\begin{eqnarray*}
 s^{-2}|(\rho_n\star\psi)(s)|^2~&=& s^{-2}
\left|\int\rho_n(s-t)^{1/2}\rho_n(s-t)^{1/2}\psi(t)~dt\right|^2~\\
&\le&
 s^{-2}\left(\frac{1}{n}+s\right)^{\!2}\int \rho_n(s-t)t^{-2}|\psi(t)|^2~dt\\
&\le&\left(1+\frac{1}{na}\right)^{\!2}~\! |(\rho_n\star(s^{-2}\psi^2))(s)|
\end{eqnarray*}
for any $s>a>0$. Since $s^{-2}\psi^2\in L^1(I_a)$ then $\rho_n\star(s^{-2}\psi^2)\to
s^{-2}\psi^2$ in $L^1(I_a)$. Thus $s^{-1}(\rho_n\star\psi)\to
s^{-1}\psi$ in $L^2(I_a)$ by the (generalized) Lebesgue Theorem, and
(\ref{eq:g1}) follows.

To prove (\ref{eq:g2}) we first argue as before to check that
$$
 s^{2}
\left|\int\rho_n(s-t)t^{-2}\psi(t)~dt\right|^2
\le \left(1-\frac{1}{na}\right)^{\!-2}~\! |(\rho_n\star(s^{-2}\psi^2))(s)|
$$
for any $s>a>0$. Thus $\rho_n\star(s^{-2}\psi)$
converges to $s^{-2}\psi$ in $L^2(I_a;s^{2}~ds)$ by Lebesgue's Theorem. In addition, $s^{-2}(\rho_n\star\psi)\to
s^{-2}\psi$ in $L^2(I_a;s^{2}~ds)$ by (\ref{eq:g1}). Thus $g_n\to 0$
in $L^2(I_a;s^{2}~ds)$ and the Lemma is completely proved. \QED

The following result for solutions to (\ref{eq:psi_intro}) is a
crucial step in the proofs of our main  theorems.

\begin{Theorem}
\label{th:ne_psi} Let $a>0$ and let $\psi$ be a distributional
solution to (\ref{eq:psi_intro}). Assume that there exists
$\gamma\le 1$ such that
$$
\psi\in L^2(I_a; s^{-2\gamma}~ds)~,\quad \alpha\ge \frac{1}{4}-(1-\gamma)^2~\!.
$$
Then $\psi=0$ almost everywhere in $I_a$.
\end{Theorem}

\proof We start by noticing that
$L^2(I_a;s^{-2\gamma}~ds)\hookrightarrow L^2(I_a; s^{-2}~ds)$ with a
continuous immersion for any $\gamma<1$. In addition, we point out
that we can assume
\begin{equation}
\label{eq:alpha}
\alpha= \frac{1}{4}-(1-\gamma)^2~\!.
\end{equation}
Let $\rho_n$ be a standard sequence of mollifiers, and let
$$\psi_n=\rho_n\star\psi~,\quad
g_n=\rho_n\star(s^{-2}\psi)-s^{-2}(\rho_n\star\psi)~\!.
$$Then
$\psi_n\to \psi$ in $L^2(I_a;s^{-2\gamma}~ds)$ and almost everywhere, and
$g_n\to 0$ in $L^2(I_a;s^2~ds)$ by Lemma \ref{L:g}. Moreover,
$\psi_n\in C^\infty(\overline I_a)$ is a non-negative solution to
\begin{equation}
\label{eq:psi_n}
-\psi''_n\ge\alpha s^{-2}\psi_n+\alpha g_n\quad\textrm{in $\mathcal D'(I_{a})$.}
\end{equation}
We assume by contradiction that $\psi\neq0$. We  let $s_0\in I_a$
such that $\eps_n:=\psi_n(s_0)\to \psi(s_0)>0$.  Up to a scaling
and after replacing $g_n$ with $s_0^2 g_n$,
we
may assume that $s_0=1$. We will show that
\begin{equation}
\label{eq:contradiction}
\eps_n:=\psi_n(1)\to \psi(1)>0
\end{equation}
leads to a contradiction.
We fix a parameter
\begin{equation}
\label{eq:delta}
\delta>\frac{1}{2}-\gamma\ge -\frac{1}{2}
\end{equation}
and for $n$ large we put
$$
\f_{\delta, n}(s):=\eps_n~\!s^{-\delta}\in L^2(I_1;s^{-2\gamma}~ds)~\!.
$$
Clearly $ \fn\in C^\infty(\R_+)$ and one easily verifies that
$(\fn)_n$ is a bounded sequence in $L^2(I_1;s^{-2\gamma}~ds)$ by  \eqref{eq:contradiction}
and (\ref{eq:delta}).
Finally we define
$$
\vn=\fn-\psi_n=\eps_n~\!s^{-\delta}-\psi_n~\!,
$$
so that $\vn\in L^2(I_1;s^{-2\gamma}~ds)$ and  $\vn(1)=0$. In addition $\vn$ solves
\begin{equation}
\label{eq:vndelta}
-\vn''\le\alpha s^{-2}~\!\vn-c_{\delta}~\!\eps_n~\!s^{-2-\delta}-\alpha g_n
\quad\textrm{in $I_1$,}
\end{equation}
where
$c_{\delta}:=\delta(\delta+1)+\alpha=\delta(\delta+1)+1/4-(1-\gamma)^2$.
Notice that $c_{\delta}>0$ and that all the
terms in the right hand side of (\ref{eq:vndelta}) belong to
$L^2(I_1;s^2~ds)$, by (\ref{eq:delta}). Thus Lemma \ref{L:D12} gives
$\vn^+\in \D(I_1)$ and
$$
\int_1^\infty|(\vn^+)'|^2~ds\le\alpha \int_1^\infty
s^{-2}|\vn^+|^2~ds-c_{\delta}\eps_n~\!\int_1^\infty
s^{-2-\delta}\vn^+~ds +o(1)~\!,
$$
since $\vn^+$ is bounded in $ L^2(I_1;s^{-2}~ds)$ and $g_n\to 0$
 in $L^2(I_1;s^2~ds)$.
 By (\ref{eq:alpha}) and  Hardy's inequality \eqref{eq:1Hardy},  we conclude that
$$  (1-\gamma)^2\int_1^\infty
s^{-2}|\vn^+|^2+c_{\delta}\,\eps_n\,\int_1^\infty
s^{-2-\delta}\vn^+~ds=o(1)~\!.
$$
Thus,  for any fixed $\delta$ we get that  $\vn^+\to 0$ almost everywhere in $I_1$  as $n\to \infty$, since $\eps_n c_\delta$ is
bounded away from $0$ by (\ref{eq:contradiction}).
 Finally we notice that
$$
\psi_n=\fn-\vn\ge \eps_n s^{-\delta}-\vn^+~\!.
$$
Since $\psi_n\to \psi$ and $\vn^+\to 0$ almost everywhere in $I_1$, and since $\eps_n\to
\psi(1)>0$, we infer that $\psi\ge \psi(1)s^{-\delta}$ in $I_1$. This conclusion
clearly
contradicts the assumption $\psi\in L^2(I_1; s^{-2\gamma}~ds)$, since
$\delta>1/2-\gamma$ was  arbitrarily chosen. Thus
(\ref{eq:contradiction}) cannot hold and the proof  is complete.
\QED

\begin{Remark}
If $\alpha>1/4$ then every non-negative solution
$\psi\in L^1_{\rm loc}(I_a)$ to problem (\ref{eq:psi_intro})
vanishes. This is an immediate consequence of  Lemma \ref{lem:AP} in
Appendix \ref{A:AP} and the sharpness of the constant $1/4$ in the
Hardy inequality (\ref{eq:1Hardy}).
\end{Remark}

\subsection{Conclusion of the proof}
\label{S:proof}
We will show that
any non-negative distributional solution $u$
to  problem
(\ref{eq:problem_ball}) gives rise to a function $\psi$ solving (\ref{eq:psi_intro}), and such that
$\psi=0$ if and only if  $u=0$. To this aim, we introduce the
Emden-Fowler transform $u\mapsto Tu$ by letting
\begin{equation}
\label{eq:EF}
 u(x) = |x|^{\frac{2-N}{2}}~(Tu)\left(\left|\log|x|\right|,\frac{x}{|x|}\right)~\!.
 \end{equation}
By change of variable formula, for any $R'\in(0,R)$ it results
\begin{equation}
\label{eq:L2} \int_{B_{R'}}\Gweight|u|^2~dx=\int_{|\log
R'|}^\infty\int_{\S^{N-1}}s^{-2\gamma}|Tu|^2~dsd\sigma~\!,
\end{equation}
 so that  $Tu\in
L^2(I_{a}\times\S^{N-1};s^{-2\gamma}~dsd\sigma)$ for any
$a>a_R:=|\log R| $. Now, for an arbitrary $\f\in C^\infty_c(I_{a_R})$
we define the radially symmetric function $\tilde\f\in
C^\infty_c(B_R)$ by setting
$$
\tilde \f(x)=|x|^{\frac{2-N}{2}}\f(|\log|x||)~\!,
$$
so that $\f=T\tilde\f$. By direct computations we get
\begin{gather}
\label{eq:CommuteT}
\int_{B_R}u(\Delta\tilde\f+\hardy|x|^{-2}\tilde\f)~dx=\int_{a_R}^\infty\f''\int_{\S^{N-1}}Tu~d\sigma ds\\
\label{eq:CommuteL}
\int_{B_R}\weight u \tilde\f~dx=\int_{a_R}^\infty s^{-2}\f \int_{\S^{N-1}}Tu~d\sigma ds~\!.
\end{gather}
Thus we are led to introduce the  function $\psi$ defined in
$I_{a_R}$ by setting
$$
 \psi(s)=
\int_{\S^{N-1}}(Tu)(s,\sigma)~\!d\sigma~\!.
$$
We notice that $\psi\in L^2(I_{a};s^{-2\gamma}~ds)$ for any $a>a_R$,
since
$$
\int_{a}^\infty s^{-2\gamma} |\psi|^2~ds\le
\left|\S^{N-1}\right|\int_{a}^\infty\int_{\S^{N-1}}s^{-2\gamma}|Tu|^2~dsd\sigma
$$
by H\"older inequality.  Moreover, from (\ref{eq:CommuteT}) and
(\ref{eq:CommuteL}) it immediately follows that $\psi\ge 0$ is a
distributional solution to
$$
-\psi''\ge \alpha s^{-2}\psi\quad\textrm{in
$\mathcal D'(I_{a_R})$.}
$$
By Theorem \ref{th:ne_psi} we infer that $\psi=0$ in $I_{a_R}$,   and
hence $u=0$ in $B_R$. The proof of Theorem \ref{th:ne_ball} is
complete. \QED

\begin{Remark}
\label{R:ce}
The assumption on  the integrability of $u$ in Theorem
\ref{th:ne_ball} are sharp. If $\alpha>1/4$ use the results in Appendix
\ref{A:AP}. For $\alpha\le 1/4$ put
$\delta_\alpha:=(\sqrt{1-4\alpha}-1)/2$ and notice that the function
$u_\alpha:B_1\to \R$ defined by
$$
u_\alpha(x)=|x|^{\frac{2-N}{2}}|\log |x||^{-{\delta_\alpha}}
$$
solves
$$
-\Delta u_\alpha -\frac{(N-2)^2}{4}|x|^{-2}~\!
u_\alpha=\alpha
|x|^{-2}\left|\log|x|\right|^{-2}~\!u_\alpha\quad \textrm{in
$\mathcal D'(B_1\setminus\{0\})$.}
$$
Moreover, if $\gamma\le 1$ then
$$
u_\alpha\in L^2_{\rm loc}(B_1;\Gweight~dx)\quad\textrm{if and only
if}\quad
\alpha<\frac{1}{4}-(1-\gamma)^2~\!.
$$
\end{Remark}

\section{Cone-like domains}
\label{S:cone}

Let $N\ge 2$. To any Lipschitz domain $\Sigma\subset \S^{N-1}$ we associate the {\em cone}
$$
\C_\Sigma:= \left\{r\sigma\in\R^N~|~\sigma\in\Sigma~,~r>0~\right\}.
$$
For any given $R>0$ we introduce also the {\em cone-like} domain
$$
\C^R_\Sigma:=\C_\Sigma\cap B_R=\left\{r\sigma\in\R^N~|~r\in(0,R)~,~\sigma\in\Sigma~\right\}.
$$
Notice that $\C_{\S^{N-1}}=\R^N\setminus\{0\}$ and
$\C_{\S^{N-1}}^R=B_R\setminus\{0\}$. If
$\Sigma$ is an half-sphere $\S^{N-1}_+$ then $\C_{\S^{N-1}_+}$ is an half-space
$\R^N_+$ and $\C_{\S^{N-1}_+}^R$ is an half-ball $B_R^+$, as in Theorem
\ref{th:half_ball}.

Assume that $\Sigma$  is properly contained in $\S^{N-1}$. Then we let $\lambda_1(\Sigma)>0$
to be the first eigenvalue of the Laplace operator on $\Sigma$.
If $\Sigma=\S^{N-1}$ we put $\lambda_1(\S^{N-1})=0$.

\noindent
It has been noticed in \cite{PT}, \cite{FaMu}, that
\begin{equation}
\label{eq:Hardy_cone}
\mu(\C_\Sigma):=\inf_
{\genfrac{}{}{0pt}{}{\scriptstyle{u\in C^\infty_c(\C_\Sigma})}{\scriptstyle{u\neq 0}}}
~\frac{\displaystyle\int_{\C_\Sigma}|\nabla u|^2~dx}
{\displaystyle\int_{\C_\Sigma}|x|^{-2}|u|^2~dx}=\frac{(N-2)^2}{4}+\lambda_1(\Sigma)~\!.
\end{equation}
The infimum $\mu(\C)$ is the best constant in the Hardy inequality
for maps having compact support in $\C_\Sigma$. In particular, for
any half-space $\R^N_+$ it holds that
$$
\mu(\R^N_+)=\frac{N^2}{4}~\!.
$$

The aim of this  section is to study the elliptic inequality
\begin{equation}
\label{eq:problemCR} -\Delta u -\mu(\C_\Sigma)|x|^{-2}~\! u\ge
\alpha |x|^{-2}\left|\log|x|\right|^{-2}~\!u\quad \textrm{in
$\mathcal D'(\CR)$.}
\end{equation}
 Notice that (\ref{eq:problemCR}) reduces
to (\ref{eq:problem_ball}) if $\Sigma=\S^{N-1}$. Problem
(\ref{eq:problemCR}) is related to an improved Hardy inequality for
maps supported in cone-like domains which will be discussed in
Appendix \ref{S:A}.

\begin{Theorem}
\label{th:ne} Let $\Sigma$ be a Lipschitz domain properly contained
in $\S^{N-1}$, $R\in(0,1]$ and  let $u\ge 0$ be a distributional
solution to (\ref{eq:problemCR}). Assume that there exists
$\gamma\le 1$ such that
$$
u\in L^2(\CR;\Gweight~dx)~,\quad \alpha\ge \frac{1}{4}-(1-\gamma)^2~\!.
$$
Then $u=0$ almost everywhere in $\CR$.
\end{Theorem}

\proof
 We introduce the first eigenfunction $\Phi\in C^2({\Sigma})\cap C(\overline{\Sigma})$ of
 the Laplace-Beltrami operator
$-\Delta_\sigma$ in $\Sigma$. Thus $\Phi$ is positive in $\Sigma$
and $\Phi$ solves
\begin{equation}
\label{eq:Phi-Sig}
\begin{cases}
-\Delta_\sigma\Phi=\lambda_1(\Sigma)\Phi&\textrm{in $\Sigma$}\\
\Phi=0~&
\textrm{on $\partial\Sigma$.}
\end{cases}
\end{equation}

Let $u\in L^2(\CR ; \Gweight~dx)$ be as in the statement, and put
$a_R=|\log R|$. We let  $Tu\in L^2(I_{a_R}\times
\Sigma;s^{-2\gamma}~dsd\sigma)$   be the Emden-Fowler transform,  as
in (\ref{eq:EF}). We further let $\psi\in
L^2(I_{a_R};s^{-2\gamma}~ds)$ defined as
$$
 \psi(s)=
\int_{\Sigma}(Tu)(s,\sigma)\Phi(\sigma)~\!d\sigma~\!.
$$
Next, for    $\f\in C^\infty_c(I_{a_R})$ being  an arbitrary
non-negative test function, we put
\begin{equation}
\label{eq:phi}
\tilde\f(x)=|x|^{\frac{2-N}{2}}\f(|\log|x|)\Phi\left(\frac{x}{|x|}\right)~\!.
\end{equation}
In essence, our aim is to test (\ref{eq:problemCR}) with $\tilde\f$
to prove that $\psi$ satisfies (\ref{eq:psi_intro}) in $I_{a_R}$. To
be  more rigorous, we use a density argument  to approximate $\Phi$
in $W^{2,2}(\Sigma)\cap H^1_0(\Sigma)$ by a sequence of smooth maps
$\Phi_n\in C^\infty_c(\Sigma)$. Then we define $\tilde\f_n$
accordingly with (\ref{eq:phi}), in such a way that
$T\tilde\f_n=\f\Phi_n$. By direct computation we get
\begin{eqnarray*}
\int_{\CR}u(\Delta\tilde\f_n+\hardy|x|^{-2}\tilde\f_n)~dx&=&
\int_{a_R}^\infty\int_{\Sigma}(Tu) \f''\Phi_n~d\sigma ds\\
&+&\int_{a_R}^\infty \int_{\Sigma}(Tu)\f\Delta_\sigma\Phi_n~d\sigma
ds
\end{eqnarray*}
\begin{gather*}
\lambda_1(\Sigma)\int_\CR |x|^{-2}u\tilde\f_n~dx=
\lambda_1(\Sigma)\int_{a_R}^\infty\int_{\Sigma}(Tu)\f\Phi_n~d\sigma ds\\
\int_{\CR}\weight u\tilde\f_n~dx= \int_{a_R}^\infty
\int_{\Sigma}s^{-2}(Tu)\f\Phi_n~d\sigma ds
\end{gather*}
Since $\tilde\f_n\in C^\infty_c(\CR)$ is an admissible test function
for (\ref{eq:problemCR}), using also (\ref{eq:Hardy_cone}) we get
\begin{eqnarray*}
-\int_{a_R}^\infty \int_{\Sigma}(Tu)\f''\Phi_n~d\sigma ds
&\ge&\alpha\int_{a_R}^\infty \int_{\Sigma}s^{-2}(Tu)\f\Phi_n~d\sigma ds\\
&-&\int_{a_R}^\infty \int_{\Sigma}(Tu)\f(\Delta_\sigma
\Phi_n+\lambda_1(\Sigma)\Phi_n)~d\sigma ds~\!.
\end{eqnarray*}
 Since $\Phi_n\to \Phi$ and
$\Delta_\sigma\Phi_n+\lambda_1(\Sigma)\Phi_n\to 0$ in
$L^2(\Sigma)$, we conclude that
$$
-\int_{a_R}^\infty \f''\psi ds
\ge\alpha\int_{a_R}^\infty s^{-2}\f\psi ds~\!.
$$
By the arbitrariness of  $\f$, we can conclude that $\psi$
 is a distributional solution to
(\ref{eq:psi_intro}). Theorem \ref{th:ne_psi} applies to give
$\psi\equiv 0$, that is, $u\equiv 0$ in $\CR$. \QED

The next result
extends Theorem \ref{th:ne} to cover the case $N=1$.
Notice that $\R_+=(0,\infty)$ is a cone and $(0,1)$ is a cone-like domain in $\R$.

\begin{Theorem}
\label{th:one}
Let $R\in(0,1]$ and let ${u}\ge 0$
 be a distributional solution to
$$
-{u}''-\frac{1}{4}~\! t^{-2}{u}\ge \alpha t^{-2}|\log
t|^{-2}~\!{u}\quad\textrm{in $\mathcal D'(0,R)$.}
$$
Assume that there exists $\gamma\le 1$ such that
$$
u\in L^2((0,R); t^{-2}|\log t|^{-2\gamma}~dt)~,\quad \alpha\ge
\frac{1}{4}-(1-\gamma)^2~\!.
$$
Then $u=0$ almost everywhere in $(0,R)$.
\end{Theorem}

\proof Write  ${u}(t)=t^{1/2}\psi\left(|\log
t|\right)=t^{1/2}\psi(s)$
 for a function
$\psi\in L^2(I_{a_R};s^{-2\gamma}~ds)$  and then notice that $\psi$
is a distributional solution to
$$
-\psi''\ge \alpha s^{-2}\psi\quad\textrm{in $\mathcal D'(I_{a_R})$.}
$$
The conclusion readily follows from Theorem \ref{th:ne_psi}.
\QED

\begin{Remark}
\label{rem:sharp_cone} If $\alpha>1/4$ then every non-negative
solution $u\in L^1_{\rm loc}(\CR)$ to problem
(\ref{eq:problemCR}) vanishes by  Theorem \ref{th:1/4cone}.

In case $\alpha\le 1/4$ the assumptions on $\alpha$
and on the integrability of $u$ in Theorems \ref{th:ne}, \ref{th:one} are sharp.
Fix $\alpha\le 1/4$, let $\delta_\alpha:=(\sqrt{1-4\alpha}-1)/2$, and
define the function
$$
u_\alpha(r\sigma)=r^{\frac{2-N}{2}}|\log r|^{-\delta_\alpha}\Phi(\sigma)~\!.
$$
Here $\Phi$ solves (\ref{eq:Phi-Sig}) if $N\ge 2$. If $N=1$ we agree
that $\sigma=1$ and $\Phi\equiv 1$. By direct computations one has
that $u_\alpha$ solves (\ref{eq:problemCR}).  Moreover, if
$\gamma\le 1$ and $R\in(0,1)$ then $u_\alpha\in
L^2(\CR;\Gweight~dx)$ if and only if
$\alpha<\frac{1}{4}-(1-\gamma)^2$.
\end{Remark}

\begin{Remark}
Nonexistence results for linear inequalities involving the differential operator
$-\Delta-\mu(\C_\sigma)|x|^{-2}$ were already obtained in \cite{FaMu}.
\end{Remark}

\appendix
\section{Hardy-Leray inequalities on cone-like domains}
\label{S:A}

In this appendix we give a simple proof of an improved Hardy
inequality for mappings having support in a cone-like domain. We
recall that for $\Sigma\subset\S^{N-1}$ we have set
$\C^1_\Sigma=\{r\sigma~|~r\in(0,1)~,~\sigma\in\Sigma~\}$ and
$\mu(\C^1_\Sigma)=(N-2)^2/4+\lambda_1(\Sigma)$.

\begin{Proposition}
\label{P:Leray}
Let $\Sigma$ be a domain in $\S^{N-1}$. Then
\begin{equation}
\label{eq:Leray}
\int_{\C^1_\Sigma}|\nabla u|^2~dx-\mu(\C_\Sigma)
\int_{\C^1_\Sigma}|x|^{-2}|u|^2\\
\ge
\frac{1}{4}\int_{\C^1_\Sigma}|x|^{-2}\left|\log|x|\right|^{-2}|u|^2~dx
\end{equation}
{for any $u\in C^\infty_c(\C^1_\Sigma)$.}
\end{Proposition}

\proof
We start by fixing an arbitrary function $v\in C^\infty_c(\R_+\times\Sigma)$.
We apply the Hardy inequality to the function
$v(\cdot,\sigma)\in C^\infty_c(\R_+)$, for any fixed $\sigma\in\Sigma$, and then we integrate over
$\Sigma$ to get
$$
\int_0^\infty\int_\Sigma|v_s|^2~dsd\sigma\ge \frac{1}{4}\int_0^\infty\int_\Sigma
s^{-2}|v|^2~dsd\sigma~\!.
$$
On the other hand, notice that $v(s,\cdot)\in C^\infty_c(\Sigma)$ for any $s\in
\R_+$. Thus, the Poincar\'e inequality for maps in $\Sigma$ plainly
implies
$$
\int_0^\infty\int_\Sigma|\nabla _\sigma v|^2~dsd\sigma- \lambda_1(\Sigma)
\int_0^\infty\int_\Sigma |v|^2~dsd\sigma\ge 0~\!.
$$
Adding these two inequalities we conclude that
$$
\int_0^\infty\int_\Sigma\left[|v_s|^2+|\nabla_{\sigma} v|^2\right]~dsd\sigma-
 \lambda_1(\Sigma)
\int_0^\infty\int_\Sigma |v|^2~dsd\sigma
\ge \frac{1}{4}\int_0^\infty\int_\Sigma
s^{-2}|v|^2~dsd\sigma
$$
for any $v\in C^\infty_c(\R_+\times\Sigma)$. We use once more the
Emden-Fowler transform $T$ in (\ref{eq:EF}) by
letting  $v:=Tu\in C^\infty_c(\R_+\times\Sigma)$ for $u\in
C^\infty_c(\mathcal C^1_\Sigma)$. Since
$$
\int_{B_1} \left[|\nabla u|^2-\hardy|x|^{-2}|u|^2\right]~\!dx=\int_0^\infty\int_{\S^{N-1}}\left[|v_s|^2+|\nabla_{\sigma}v|^2\right]~dsd\sigma~\!,
$$
then (\ref{eq:L2}) readily leads to the conclusion.
\QED

\begin{Remark}\label{rem:sharp-Cone}
The arguments we have used to prove Proposition \ref{P:Leray} and
the fact that the best constant in the Hardy inequality for maps in
$C^\infty_c(\R_+)$ is not achieved show that the constants in
inequality (\ref{eq:Leray}) are sharp, and not achieved.
\end{Remark}

\begin{Remark}\label{rem:sharp-B}
\label{R:HL} Notice that for $N\geq1$, we have $\mathcal
C_{\S^{N-1}}=\R^N\setminus\{0\}$ and $\mu(\mathcal
C_{\S^{N-1}})=(N-2)^2/4$ . Thus \ref{eq:Leray} gives  (\ref{eq:HL})
for $u\in C^\infty_c(B_1\setminus\{0\})$.
\end{Remark}

In the next proposition we extend the inequality (\ref{eq:Leray}) to cover the case
$N=1$.

\begin{Proposition}
\label{P:Leray1}
It holds that
$$
\int_0^1|{u}'|^2~dt-\frac{1}{4}\int_0^1t^{-2}|{u}|^2~dt\ge
\frac{1}{4} \int_0^1t^{-2}|\log t|^{-2}|{u}|^2~dt
$$
for any ${u}\in C^\infty_c(0,1)$. The constants are sharp, and not achieved.
\end{Proposition}

\proof Write  ${u}(t)=t^{1/2}\psi\left(|\log
t|\right)=t^{1/2}\psi\left(s\right)$ for a function $\psi\in
C^\infty_c(\R_+)$ and then apply the Hardy inequality to $\psi$.
\QED

Next, let $\theta\in\R$ be a given parameter and let $\Sigma$ be a Lipschitz
domain in $\S^{N-1}$, with $N\ge 2$. For an arbitrary
$u\in C^\infty_c(\mathcal C_\Sigma^1)$ we put $v=|x|^{-\theta/2}u$. Then the Hardy-Leray inequality
(\ref{eq:Leray}) and integration by parts plainly imply that
$$
\int_{\C^1_\Sigma}|x|^{\theta}|\nabla v|^2~dx-\mu(\C_\Sigma;\theta)
\int_{\C^1_\Sigma}|x|^{\theta-2}|v|^2\\
\ge
\frac{1}{4}\int_{\C^1_\Sigma}|x|^{\theta-2}\left|\log|x|\right|^{-2}|v|^2~dx
$$
for any $v\in C^\infty_c(\C^1_\Sigma)$, where
\begin{equation}
\label{R:a} \mu(\C_\Sigma;
\theta):=\frac{(N-2+\theta)^2}{4}+\lambda_1(\Sigma)~\!.
\end{equation}
It is well known that
$$
\frac{(N-2+\theta)^2}{4}=
\inf_
{\genfrac{}{}{0pt}{}{\scriptstyle{u\in C^\infty_c(\R^N\setminus\{0\})}}{\scriptstyle{u\neq 0}}}
~\frac{\displaystyle\int_{B_1}|x|^\theta|\nabla u|^2~dx}
{\displaystyle\int_{B_1}|x|^{\theta-2}|u|^2~dx}
$$
is the Hardy constant relative to the operator $L_\theta v=-\div(|x|^\theta\nabla v)$.
For the case $N=1$ one can obtain in a similar way the inequality
$$
\int_0^1t^\theta|v'|^2~dt-\frac{(\theta-1)^2}{4}\int_0^1
t^{\theta-2}|v|^2~dt\ge \frac{1}{4}\int_0^1 t^{\theta-2}|\log
t|^{-2}|v|^2~dt~\!
$$
which holds for any $\theta\in\R$ and for any $v\in
C^\infty_c(0,1)$.

\section{A general necessary condition}
\label{A:AP} In this appendix we show in particular that a necessary
condition for the existence of non-trivial and non-negative
solutions to  (\ref{eq:problem_ball}) and (\ref{eq:problemCR}) is
that $\alpha\le 1/4$. We need the following general lemma, which
naturally fits into the classical Allegretto-Piepenbrink theory (see
for instance \cite{Al} and \cite{Pi}).

\begin{Lemma}\label{lem:AP}
Let $\O$ be a domain in $\R^N$, $N\geq1$. Let $a\in
L^\infty_{loc}(\O)$ and $a>0$ in $\O$. Assume that $u\in
L^1_{loc}(\O)$ is a non-negative, non-trivial solution to
$$
-\Delta u\geq a(x) u\quad \mathcal D'(\O).
$$
Then
$$
\int_{\O}|\n \phi|^2~dx\geq
\int_{\O}a(x)\,|\phi|^2~dx,\quad\textit{for any $\phi\in
C^\infty_c(\O)$.}
$$
\end{Lemma}

\proof
Let  $A\subset\O$ be a measurable set such that $|A|>0$ and $u>0$ in $A$.
Fix any function
$\phi\in C^\infty_c({\O})$ and choose a domain $\wt\O\subset\subset\O$ such that $|\wt{\O}\cap A|>0$
and $\phi\in C^\infty_c({\wt\O})$.
For any integer $k$ large enough put $f_k=\min\{a(x)u,k\}\in
L^\infty(\wt{\O})$. Let $v_k\in H^1_0(\wt\O)$ be the unique solution to
\begin{equation}
\label{eq:problem_v}
\begin{cases}
-\Delta v_k=f_k&\textrm{in $\wt{\O},$}\\
v_k=0&\textrm{on $\de \wt{\O}$.}
\end{cases}
\end{equation}
Notice that $v\in C^{1,\b}(\wt{\O})$ for any $\b\in(0,1)$. Since for $k$ large enough the function $f_k$ is non-negative
and non-trivial then $v\ge 0$. Actually it turns out that $v^{-1}\in L^\infty_{loc}(\wt{\O})$
by the Harnack inequality. Finally, a convolution argument and the
maximum principle plainly give
\begin{equation}\label{eq:mmp}
u\ge v_k>0\quad   \textrm{almost everywhere in $\wt{\O}$}.
\end{equation}
Since $v_k^{-1}\phi\in L^\infty(\wt{\O})$ then we can use
$v_k^{-1}\phi^2$ as test function for (\ref{eq:problem_v}) to get
$$
\int_\O\nabla v_k\cdot\nabla\left(v_k^{-1}\phi^2\right)~dx=\int_\O f_k v_k^{-1}\phi^2~dx
\ge \int_\O f_k u^{-1}\phi^2~dx
$$
by \eqref{eq:mmp}. Since
$\nabla v_k\cdot\nabla\left(v_k^{-1}\phi^2\right)
=|\nabla \phi|^2-\left|v_k\nabla(v_k^{-1}\phi)\right|^2\le|\nabla\phi|^2$,
we readily infer
$$
 \int_{\O}|\n \phi|^2~dx\geq  \int_{\O}f_k u^{-1} \phi^2~dx
$$
and Fatou's lemma implies that
$$
\int_{\O}|\nabla \phi |^2~dx \ge \int_{\O}a(x)\,\phi^2~dx~\!.
$$
The conclusion readily follows.
 \QED

The sharpness of the constants in (\ref{eq:HL}) (compare with Remark \ref{rem:sharp-Cone})
and  Lemma \ref{lem:AP} plainly imply the following
result.

\begin{Theorem}
\label{th:1/4}
Let $N\geq1$, $R\in(0,1]$ and $c, \alpha\geq0$. Let $u\in
L^1_{\textrm loc}(B_R\setminus\{0\})$ be a non-negative
distributional solution to
$$
-\Delta u -c~\!|x|^{-2}~\! u\ge \alpha |x|^{-2}\left|\log|x|\right|^{-2}~\!u\quad
\textrm{in $\mathcal D'(B_R\setminus\{0\})$.}
$$

$i)$~ If $c>\hardy$ then $u\equiv0$.

$ii)$ If $c=\hardy$ and $\alpha>\frac{1}{4}$ then $u\equiv0$.
\end{Theorem}
We notice that proposition $i)$ in Theorem \ref{th:1/4} was already proved in
\cite{BG} (see also  \cite{Du}).

Finally, from
Remark \ref{rem:sharp-Cone} and Lemma \ref{lem:AP}, we obtain the next nonexistence result.
\begin{Theorem}
\label{th:1/4cone}
Let $\Sigma$ be a domain properly contained in $\S^{N-1}$,
$R\in(0,1]$ and $c, \alpha\geq0$. Let $u\in L^1_{\textrm loc}(\CR)$
be a non-negative distributional solution to
$$
-\Delta u -c~\!|x|^{-2}~\! u\ge \alpha |x|^{-2}\left|\log|x|\right|^{-2}~\!u\quad
\textrm{in $\mathcal D'(\CR)$.}
$$

$i)$~ If $c>\mu(\calC_\Sig)$ then $u\equiv0$.

$ii)$ If $c=\mu(\calC_\Sig)$ and $\alpha>\frac{1}{4}$ then $u\equiv0$.
\end{Theorem}

\section{Extensions}
\label{S:B}

In this appendix we state some nonexistence theorems that can proved
by using a suitable functional change $u\mapsto \psi$ and Theorem
\ref{th:ne_psi}. We shall also point out some corollaries of our
main results.

\subsection{The $k$-improved weights}
We define a sequence of radii $R_k\to0$ by setting $R_1=1$,
$R_k=e^{-\frac{1}{R_{k-1}}}$. Then we use induction again to define
 two sequences of radially symmetric weights $X_k(x)\equiv X_k(|x|)$ and $z_k$ in $B_{R_k}$ by setting
$X_1(|x|)=|\log |x||^{-1}$ for $|x|<1=R_1$ and
$$
X_{k+1}(|x|)=X_k\left(|\log |x||^{-1}\right)~,\quad z_k(x)=|x|^{-1}\prod_{i=1}^{k}X_i(|x|)
$$
for all $x\in B_{R_k}\setminus \{0\}$. It can be proved by induction that $z_k$ is well defined on
$B_{R_k}$ and $z_k\in L^2_{\rm loc}(B_{R_k})$.
We are interested in distributional solutions to
\begin{equation}
\label{eq:problem_k} -\Delta u-\hardy |x|^{-2}u\geq \alpha
z_k^2~\! u\quad \mathcal{D}'(B_R\setminus\{0\})~\!
\end{equation}
for $R\in(0,R_k]$. The next result includes Theorem \ref{th:ne_ball}
by taking $k=1$.

\begin{Theorem}
\label{th:ne_k}
Let $k\ge 1$, $R\in(0,R_k]$ and let $u\ge 0$ be a distributional solution to
(\ref{eq:problem_k}). Assume that there exists $\gamma\le 1$ such that
$$
u\in L^2_{\rm loc}(B_R;z_k^2X_k^{2(\gamma-1)}~dx)~,\quad
\alpha\ge \frac{1}{4}-(1-\gamma)^2~\!.
$$
Then $u=0$ almost everywhere in $B_R$.
\end{Theorem}

\proof
We start by introducing the $k^{\it th}$ Emden-Fowler transform $u\mapsto T_ku$,
$$
u(x)=z_{k}(|x|)^{-\frac{1}{2}}|x|^{\frac{1-N}{2}}X_k(|x|)^{\frac{1}{2}}
(T_ku)\left(X_k(|x|)^{-1},\frac{x}{|x|}\right)~\!.
$$
Notice that for any $R<R_k$ it results \be \label{eq:zk_integral}
\int_{B_{R}} z_{k}^2~\!X_k^{2(\gamma-1)} |u|^{2}~dx=
\int_{X_k(R)^{-1}}^\infty
s^{-2\gamma}\int_{\S^{N-1}}|T_ku|^2~dsd\sigma~\!, \ee
so that $T_ku\in L^2(I_a\times\S^{N-1}; s^{-2\gamma}~dsd\sigma)$ for any $a>X_k(R)^{-1}$.
This can be easily checked by noticing that $X'_k=z_k X_k$.
Next we set
$$
\psi_u(s):=\int_{\S^{N-1}}(T_ku)(s,\sigma)d \sigma~\!.
$$
By (\ref{eq:zk_integral}) we have that  $\psi\in
L^2(I_a;s^{-2\gamma}~ds)$ for any $a>X_k(R)^{-1}$. Thanks to Theorem
\ref{th:ne_psi}, to conclude the proof it suffices to show that
$\psi$ is a distributional solution to $-\psi''\ge\alpha s^{-2}\psi$
in the interval $I_{\tilde a}$, where $\tilde a= X_k(R)^{-1}$. To
this end, fix any test function $\f\in C^\infty(I_{\tilde a})$, and
define the radially symmetric mapping $\tf\in
C^\infty_c(B_R\setminus\{0\})$ such that $T_k\tf=\f$. By direct
computation one can prove that
$$
\Delta\tf+\hardy|x|^{-2}\tf=\omega\tf+|x|^{\frac{1-N}{2}}z_k^{\frac{3}{2}}X_k^{-\frac{3}{2}}
\f''\left(X_k(|x|)^{-1}\right)
$$
where $\omega\equiv 0$ if $k=1$, and
$$
\omega=\frac{1}{2}\left[
\left(\sum_{i=1}^{k-1}z_i\right)^2-\frac{1}{2}\sum_{i=1}^{k-1}z_i^2
 \right]~\!
$$
if $k\ge 2$. Since $\omega\ge 0$ then
$$
\int_{B_R}u\left(\Delta\tf+\hardy|x|^{-2}\tf\right)~dx\ge\int_{\tilde a}^\infty \psi
\f''~\!ds
$$
provided that $\f$ is non-negative. In addition it results
$$
\int_{B_R}z_k^2~\!u\tf~dx=\int_{\tilde a}^\infty s^{-2}\psi\f~ds~\!.
$$
Since $\f$ was arbitrarily chosen, the conclusion readily
follows. \QED

By similar arguments as above and in Section~\ref{S:cone}, we can
prove a nonexistence result of positive solutions to the  problem
\begin{equation}
\label{eq:problem_k-cone} -\Delta u-\mu(\C_\Sigma) |x|^{-2}u\geq \alpha
z_k^2~\! u\quad \mathcal{D}'(\CR)~\!,
\end{equation}
where $\C_\Sigma$ is a Lipschitz proper cone in $\R^N$, $N\ge 1$, and
$C_\Sigma^R=\C_\Sigma\cap B_R$. We shall skip the proof the following result.

\begin{Theorem}
\label{th:ne_k-cone} Let $k\ge 1$, $R\in(0,R_k]$ and let $u\ge 0$ be
a distributional solution to \eqref{eq:problem_k-cone}. Assume that
there exists $\gamma\le 1$ such that
$$
u\in L^2_{ }(\CR;z_{k}^2X_k^{2(\gamma-1)}~dx)~,\quad
\alpha\ge \frac{1}{4}-(1-\gamma)^2~\!.
$$
Then $u=0$ almost everywhere in $\CR$.
\end{Theorem}
Some related improved Hardy inequalities involving the weight $z_k$
and which motivate the interest of problems (\ref{eq:problem_k}) and
(\ref{eq:problem_k-cone}) can be found in  \cite{ACR}, \cite{Ch},
\cite{GM} and also \cite{BFT}.

\subsection{Exterior cone-like domains}
The Kelvin transform
$$
u(x)\mapsto |x|^{2-N}u\left(\frac{x}{|x|^2}\right)
$$
can be used to get nonexistence results for
exterior domains in $\R^N$.

Let $\Sigma$ be a domain in $\S^{N-1}$, $N\ge 2$, and let $\C_\Sigma$ be
the cone defined in Section \ref{S:cone}. We recall that
$\mu(\C_\Sigma)=(N-2)^2/4+\lambda_1(\Sigma)$. Since the inequality in
(\ref{eq:problem_ball}) is invariant with respect to the Kelvin transform, then
Theorems \ref{th:ne_ball} and \ref{th:ne} readily lead to the following
nonexistence result.

\begin{Theorem}
\label{th:ne_ext}
Let $\Sigma$ be a Lipschitz domain in $\S^{N-1}$, with $N\ge 2$. Let
$R>1$,  $\alpha\in\R$ and let $u\ge 0$ be a distributional solution to
$$
-\Delta u -\mu(\C_\Sigma)|x|^{-2}~\! u\ge \alpha |x|^{-2}\left|\log|x|\right|^{-2}~\!u\quad
\textrm{in $\mathcal D'(\C_\Sigma\setminus \overline B_R)$.}
$$
Assume that there exists $\gamma\le 1$ such that
$$
u\in L^2(\C_\Sigma\setminus \overline B_R;|x|^{-2}|\log|x||^{-2\gamma}~dx)~,
\quad \alpha\ge \frac{1}{4}-(1-\gamma)^2~\!.
$$
Then $u=0$ almost everywhere in $\C_\Sigma\setminus \overline B_R$.
\end{Theorem}

A similar statement holds in case $N=1$ for ordinary differential
inequalities in unbounded intervals $(R,0)$ with $R>0$, and for
problems involving the weight $z_k^2$.

\subsection{Degenerate elliptic operators}

Let $\theta\in \R$ be a given real parameter. We notice that $u$ is
a distributional solution to (\ref{eq:problemCR}) if and only if
$v=|x|^{-\theta/2}u$ is a distributional solution to
\begin{equation}
\label{eq:problem_a} -\div(|x|^ \theta\nabla
v)-\mu(\C_\Sigma;\theta) |x|^{\theta-2}~\!v\ge \frac{1}{4}~\!
|x|^{\theta-2}|\log|x||^{-2}~\!v\quad\textrm{in $\mathcal
D'(\CR)$}~,
\end{equation}
where   $\mu(\C_\Sigma;\theta)$ is defined in Remark \ref{R:a}.
Therefore Theorem \ref{th:ne_ball} and Theorem \ref{th:ne} imply the
following nonexistence result for linear inequalities involving the
weighted Laplace operator $L_\theta v=-\div(|x|^\theta\nabla v)$.

\begin{Theorem}
\label{th:a} Let $\Sigma$ be a Lipschitz domain in $\S^{N-1}$. Let
$\theta\in\R$, $R\in(0,1]$,  $\alpha\in\R$ and let $v\ge 0$ be a
distributional solution to (\ref{eq:problem_a}). Assume that there
exists $\gamma\le 1$ such that
$$
v\in L^2(\CR;|x|^{\theta-2}|\log|x||^{-2\gamma}~dx)~, \quad
\alpha\ge \frac{1}{4}-(1-\gamma)^2~\!.
$$
Then $v=0$ almost everywhere in $\CR$.
\end{Theorem}

A nonexistence result for the operator $-\div(|x|^\theta\nabla v)$
similar to Theorem \ref{th:ne_ext} or to Theorem \ref{th:ne_k}
 can be obtained from Theorem \ref{th:a}, via  suitable functional
 changes.

 \label{References}

\end{document}